\documentclass[10pt,journal]{IEEEtran}

\usepackage{graphicx}
\usepackage[cmex10]{amsmath}
\usepackage{amsfonts} 
\usepackage{amssymb}
\usepackage{stmaryrd}
\usepackage{subfigure}

\def\R{{\mathbb{R}}}

\def\Div{\textup{div\,}}

\newtheorem{thm}{Theorem}
\newtheorem{lem}{Lemma}

\begin{document}
\title{Properties of $BV-G$ structures+textures decomposition models. Application to road detection in satellite images}

\author{J\'er\^ome~Gilles and~Yves~Meyer,
\thanks{J. Gilles is with EORD Department of DGA/CEP, 16bis rue Prieur de la C\^ote d'Or, 94110 Arcueil, France, email: jerome.gilles@etca.fr (see http://jerome.gilles91.free.fr).}%
\thanks{Y. Meyer is with CMLA-ENS Cachan, 62 av du Pr\'esident Wilson 94230 Cachan, France, email: yves.meyer@cmla.ens-cachan.fr.}%
\thanks{Manuscript received June xx, 2009; revised xx.}}

\markboth{Journal of Trans on Image Processing,~Vol.~6, No.~1, January~2010}%
{\MakeLowercase{Theoretical properties of $BV-G$ structures+textures decomposition models}}%

\maketitle

\begin{abstract}
In this paper we present some theoretical results about a structures-textures image decomposition model which was proposed by the second author. We prove a theorem which gives the behavior of this model in different cases. Finally, as a consequence of the theorem we derive an algorithm for the detection of long and thin objects applied to a road networks detection application in aerial or satellite images.
\end{abstract}



\section{Introduction}
A few years ago, the second author proposed several image decomposition models which are aimed to split an image into three components. The first component should describe the objects contained in the image, the second one is given by the textured components and the third one is an additive noise. The structures are modelized as functions belonging to the $BV-$space (the space of bounded variations functions) and the textures as oscillating functions belonging to a space, called $G$ by the second author, which is close to the dual of the $BV-$space. Many papers deal with numerical implementation \cite{aujol,vese2,vese1}, extensions to other cases (like color images \cite{aujolcouleur}) and some applications \cite{aujolclassiftexture} but few of them study the theoretical behavior of this model.\\

In this paper, we propose to explore this behavior through different results. The main theorem we present depicts the optimal decompositions obtained by an adapted tuning of parameters and the properties in the different functional spaces of the image to be decomposed.\\

The remainder of the paper is organized as follows. In Section \ref{sec:tvl2}, we remind the Rudin-Osher-Fatemi algorithm which is the origin of the work of the second author about image decomposition model. In Section \ref{sec:tvg}, we present the $BV-G$ model and give some notations which will be used in the rest of the paper. Section \ref{sec:new} is the main section of the paper and presents some theoretical results detailing the behavior of the model. In Section \ref{sec:app}, we present an application of the main theorem proved in Section \ref{sec:new} to the enhancement of long and thin structures which can be used for example for road networks detection in aerial or satellite images. We will finish by concluding and giving some perspectives to this work for future research.

\section{$BV-L2$ model - The Rudin Osher Fatemi algorithm}\label{sec:tvl2}
The starting point of the second author's work is the Rudin-Osher-Fatemi (ROF) model \cite{rof}. At the origin, this algorithm was developped for image restoration purposes. It gives good results and is nowadays currently used. The authors propose to retrieve the restored image $u$ from the corrupted image $f$ by assuming that $u$ belongs to the space $BV$, the space of bounded variations functions which is well adapted to modelize structures in an image and widely used in the literature. They propose to minimize the functional (\ref{eq:rof}).

\begin{equation}\label{eq:rof}
F_{\lambda}^{ROF}(u)=J(u)+\lambda \|f-u\|_{L^2}^2
\end{equation}
where $J(u)=\int |\nabla u|$ is the total variation (TV) of $u$ (this corresponds to the fact we want $u$ belonging to $BV$) and $\lambda$ is a regularization parameter. \\
In \cite{chambolle} the author proposes, in the case of a numerical framework and for bounded domain, a very efficient nonlinear projector (denoted $P_{G_{\lambda}}$) to find the minimizer $\hat{u}$ of (\ref{eq:rof}). Then $\hat{u}=f-P_{G_{\lambda}}(u)$ and the advantages of this algorithm are twofold:
\begin{enumerate}
\item it is very easy to implement by an iterative process,
\item a theorem gives the condition which ensures the convergence in the discrete case.\\
\end{enumerate}

Now, we take the point of view of image decomposition by assuming that $f=u+v$ where $v$ corresponds to the textures of the image. We mean by textures some oscillating patterns in the image. Then we can rewrite the ROF model as (\ref{eq:rofbis}) (in practice, the minimization occurs only on $u$ and $v$ is obtained by $v=f-u$).
\begin{equation}\label{eq:rofbis}
F_{\lambda}^{ROF}(u)=J(u)+\lambda \|v\|_{L^2}^2.
\end{equation}
But this model is not adapted to capture textures. For example, let us assume that $g$ is a texture defined by $g(x)=\theta(x)\cos(Nx_1)$, where $\theta(x)$ is the characteristic function on the unit square, $N$ is the frequency and $x_1$ a direction in the image. Then, we can check that 
\begin{equation}
\|g\|_{L^2}\approx \frac{1}{\sqrt{2}}\|\theta\|_{L^2}
\end{equation}
which does not depend on $N$, and
\begin{equation}
\|g\|_{BV}=\frac{2N}{\pi}\|\theta\|_{L^1}+\epsilon_N
\end{equation}
which tends to infinity when $N\rightarrow\infty$ ($\epsilon_N\rightarrow 0$ when $N\rightarrow\infty$). So more some patterns are oscillating, the less the algorithm is relevant to correctly capture them.\\

More generally, the ROF model has three main defects.\\

First, in the continuous case the model has no meaning if the image is corrupted by white noise. Indeed, the $L^2$-norm of a white noise is infinite. In dimension two, a gaussian white noise has a finite norm only in function spaces with negative regularity indices.

The second one is that textures and noise are treated in the same manner whereas textures is somewhat ``structured'' (like periodicity or high frequencies for example) and noise is completely unstructured. In term of Fourier analysis, some localized frequencies exist for textures while gaussian white noise has a constant spectrum.

The last one, as shown in \cite{meyer}, it exists a slight loss of intensity when $f$ is a constant times the characteristic function of a disk.

Based on these remarks, the second author proposed in \cite{meyer} to use some specific functional spaces and their associated norms to model textures. This approach will be described in the next section.

\section{$BV-G$ decomposition model}\label{sec:tvg}
In this section we present the model proposed by the second author \cite{meyer} to decompose an image $f$ into two parts, structures $u$ and textures $v$. We saw in the previous section that the Rudin-Osher-Fatemi model is not adapted to deal with the texture component. The author proposes to modify the ROF functional as in (\ref{eq:meyer}).
\begin{equation}\label{eq:meyer}
F^{YM}(u,v)=J(u)+(2\lambda)^{-1}\|v\|_G,
\end{equation}

In order to introduce the definition of the space $G$, let us recall some properties of the space $BV$. It exists different (equivalent) ways to define the space $BV$ (see \cite{chambolle2} for an introduction on total variation). The total variation can be defined by duality: $\forall u\in L_{loc}^1(\Omega)$, the total variation is given by (\ref{eq:deftv}).
\begin{align}\label{eq:deftv}
J(u)=\sup\bigg\{-\int_{\Omega} u\Div\phi dx: &\phi\in\mathcal{C}_c^{\infty}(\Omega,\R^N), \\ \notag
& |\phi|\leqslant 1\; \forall x\in\Omega \vphantom{\int} \bigg\}
\end{align}
which, $\forall u\in \mathcal{C}^1$, is equivalent to $J(u)=\int |\nabla u|$. Then the space of bounded variation functions, $BV$, is endowed by the norm $\|.\|_{BV}=\|.\|_{L^1}+J(.)$. The dual of $BV$ is not a functional space, but if we consider the closure of $\mathcal{S}(\R^2)$ in $BV$ (which is denoted $\mathcal{BV}$), the dual of $\mathcal{BV}$ is a functional space denoted $G$. The $G-norm$ is defined by the following recipe.
For $v=\Div g=\partial_1 g_1+\partial_2 g_2$ where $g=(g_1,g_2)\in L^{\infty}(\mathbb{R}^2)\times L^{\infty}(\mathbb{R}^2)$,
\begin{equation}\label{eq:normeg}
\|v\|_G=\inf_{g}\left\|\left(\left|g_1\right|^2+\left|g_2\right|^2\right)^{\frac{1}{2}}\right\|_{L^{\infty}}.
\end{equation}

Even if there's no direct duality relation between $BV$ and $G$, it is easy to see that the $G-norm$ and the total variation have dual behaviors. We mean that $BV$ is devoted to modelize structures in an image (like characteristic functions for example) and the space $G$ is well-adapted to modelize oscillating patterns. We saw in the previous section that for $g(x)=\theta(x)\cos(Nx_1)$, the $BV-norm$ is not adapted to capture it. We can easily check that $\|g\|_G \leqslant \frac{C}{N}$ which confirms that the oscillating pattern will be captured by the $G-norm$.\\ Now, if we return to the functional (\ref{eq:meyer}) proposed by the second author, we can easily understand its behavior. Assume we want to minimize this functional over $u$ and $v$, the $BV-norm$ reaches its minimum if $u$ corresponds to structures in the image and the $G-norm$ reaches its minimum for oscillating patterns (like textures). Then this model permits a better structures + textures decomposition than the ROF model presented in the previous section. \\

\begin{figure*}[!t]
\centerline{
\subfigure[Original]{\includegraphics[scale=0.33]{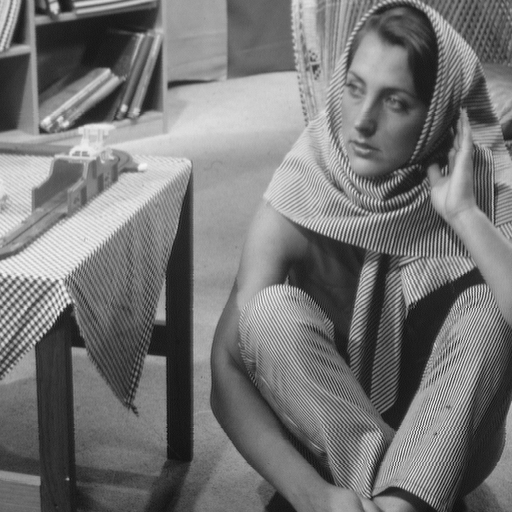}}
\hfill
\subfigure[Structures]{\includegraphics[scale=0.33]{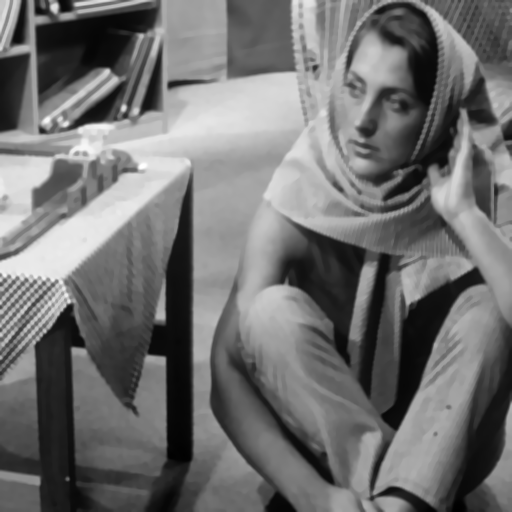}}
\hfill
\subfigure[Textures]{\includegraphics[scale=0.33]{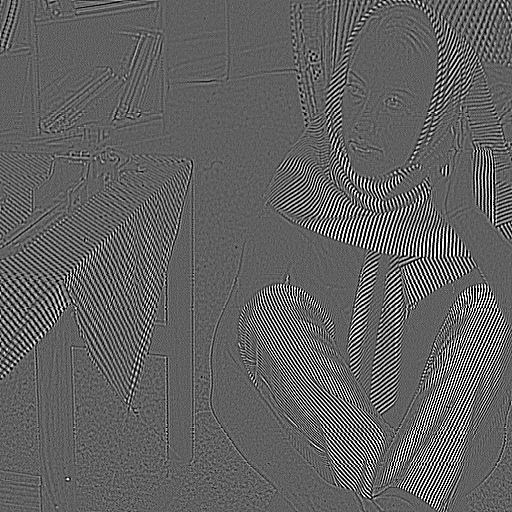}}
}
\caption{Example of image structures+textures decomposition.}
\label{fig:uv}
\end{figure*}

Vese and Osher, in \cite{vese1} were the first ones to propose numerical experiments of the second author's model. In \cite{aujol,aujolphd}, Aujol proposes to use the nonlinear projector of Chambolle to solve the model. It consists of an alternate iterative algorithm which provides the minimizers of (\ref{eq:meyer}). These minimizers are given by $\hat{u}=f-\hat{v}-P_{G_{\lambda}}(f-\hat{v})$ and $\hat{v}=P_{G_{\mu}}(f-\hat{u})$ where $\mu$ is an upper bound for $\|v\|_G$. Figure \ref{fig:uv} shows one example of a structure+texture decomposition of an image.\\
In \cite{aujoluvw,jegilles,triet1,triet2} the authors propose different approaches to extend this model to a three components model to deal with noisy images. These models permit to separate the structures, textures and noise respectively. In this paper we restrict our study to the two components model.

We begin the next section by giving new theoretical results.

\section{Decomposition model output}\label{sec:new}

\subsection{General results}
In this subsection we first recall some important results given by the second author (proofs can be found in \cite{meyer}) which will be useful to our work. The first one concerns the following inequality

\begin{lem}\label{lemme:bv1}
If $u\in L^2(\R^2)$ and $v\in BV(\R^2)$, then
\begin{equation}
\left|\int u(x)v(x)dx\right| \leqslant \|u\|_G \|v\|_{BV}.
\end{equation}
\end{lem}

The next two lemmas are general results on functional minimization in Banach$/L^2$-spaces and Banach$/$Banach-spaces respectively.
\begin{lem}\label{lemme:meyer1}
Let $E$ be an arbitrary Banach space and let $\|.\|_{E^*}$ be the dual norm. Let us assume we want to minimize
\begin{equation}
\|u\|_E+\lambda\|v\|_{L^2}^2
\end{equation}
for all decompositions $f=u+v$ of $f\in L^2(\R^2)$, then two cases appear\\
\makeatletter
\renewcommand{\labelenumi}{(\theenumi)}
\makeatother
\begin{enumerate}
\item if $\|f\|_{E^*}\leqslant \frac{1}{2\lambda}$, then the minimum is reached for $u=0$ and $v=f$,\label{enumaa} \\
\item if $\|f\|_{E^*}> \frac{1}{2\lambda}$, then the minimum is reached for $v$ such that $\|v\|_{E^*}=\frac{1}{2\lambda}$ and $\langle u,v \rangle=\frac{1}{2\lambda}\|u\|_E$.\label{enumbb} \\
\end{enumerate}
\end{lem}

\begin{lem}\label{lemme:meyer2}
Let $E_1$ and $E_2$ two Banach spaces embedded in a vector space $E$. We define the Banach space $E_3$ which is the set of all $z$ such that 
\begin{equation}
z=x+y\; , \; x\in E_1 \; , \; y\in E_2
\end{equation}
provided with the norm
\begin{equation}
\|z\|_{E_3}=\inf \{\|x\|_{E_1}+\|y\|_{E_2}\}.
\end{equation}
Then $E_3$ is the smallest Banach space containing $E_1$ and $E_2$. In addition, $E_3^*$ is the biggest Banach space contained in $E_1^*$ and $E_2^*$. In other terms, $E_3^*=E_1^*\cap E_2^*$ and the norm of $g$ in $E_3^*$ is defined by
\begin{equation}
\|g\|_{E_3^*}=\sup \{\|g\|_{E_1^*}\; , \; \|g\|_{E_2^*}\}.
\end{equation}
\end{lem}

\subsection{Decomposition model properties}
The space $G$ is defined as above (Section \ref{sec:tvg}). To study the $BV-G$ decomposition model, we propose to use the following approach. Let us assume that we are given an image $f\in L^2(\R^2)$ and two positive parameters $\lambda$ et $\mu$. Then we seek to decompose $f$ as a sum
\begin{equation}
f=u+v+w
\end{equation}
by minimizing the functional $E(u,v,w)$ defined by 
\begin{equation}\label{eq:pbmeyer}
\|u\|_{BV}+\lambda \|v\|_{L^2}^2+\mu\|w\|_G.
\end{equation}

The ROF model corresponds to the case $\mu=+\infty$. As $BV\subset L^2$, we necessarily have $w\in L^2$. 
The existence of an optimal decomposition is given by the ``Hilbert's direct method''. Since $BV$ is a dual space, from every bounded sequence $u_j\in BV$ we can extract a subsequence that converges, in the distributionnal sense, to $u\in BV$. The same argument can be used for $L^2$ and $G$. The uniqueness is not ensured except for the $v$ part. More details will be given later in this paper.\\

Before demonstrating the main theorem, we introduce some intermediate results.

\begin{lem}\label{lem:isoperi}
For all function $f\in BV$, we have
\begin{equation}
\|f\|_{L^2}\leqslant \frac{1}{2\sqrt{\pi}}\|f\|_{BV}
\end{equation} 
and this implies 
\begin{equation}
\|f\|_G\leqslant \frac{1}{2\sqrt{\pi}}\|f\|_{L^2}\leqslant \frac{1}{4\pi}\|f\|_{BV}.
\end{equation}
\end{lem}

This lemma is a direct consequence of the isoperimetric inequality.

We deduce the following theorem.\\

\begin{thm}
\label{theo:meyer2}
If $0<\mu<4\pi$, then the optimal decomposition $f=u+v+w$ verifies $u=0$.\\
\end{thm}

\begin{IEEEproof}
Assume that we fix $v$ and let $u$ free. We will write $u+w=\sigma$ and then $\sigma=f-v$. First, we want to minimize $\|u\|_{BV}+\mu\|\sigma-u\|_G$. If we assume that $0<\mu<4\pi$, by Lemma \ref{lem:isoperi} we get

\begin{align}
\|u\|_{BV}+\mu \|\sigma-u\|_G & \geqslant 4\pi\|u\|_G+\mu\|\sigma-u\|_G\\
& \geqslant \mu\|u\|_G+\mu\|\sigma-u\|_G \\
& \geqslant \mu\|\sigma\|_G.
\end{align}
In addition, if $u$ is not zero then $\|u\|_G>0$ and we have
\begin{equation}
\|u\|_{BV}+\mu\|\sigma-u\|_G>\mu\|\sigma\|_G.
\end{equation}
We conclude that the minimum will be reached for $u=0$.
\end{IEEEproof}

The following theorem gives the behavior of the model described by (\ref{eq:pbmeyer}).\\

\begin{thm}\label{theo:meyer1}
If $\|f\|_G\leqslant \frac{1}{2\lambda}$ and $\|f\|_{BV}\leqslant \frac{\mu}{2\lambda}$, then $u=w=0$ and the optimal decomposition is $f=0+f+0$.\\

If $\|f\|_G\leqslant \frac{1}{2\lambda}$ but $\|f\|_{BV}>\frac{\mu}{2\lambda}$, three cases appear for an optimal decomposition $f=u+v+w$. \\

\makeatletter
\renewcommand{\labelenumi}{(\theenumi)}
\makeatother
\begin{enumerate}
\item $u=0$, $\|v\|_{BV}=\frac{\mu}{2\lambda}$, $\|v\|_G<\frac{1}{2\lambda}$ and $\langle v,w \rangle=\frac{\mu}{2\lambda}\|w\|_G$, \label{enuma} \\
\item $w=0$, $\|v\|_{BV}\leqslant\frac{\mu}{2\lambda}$, $\|v\|_G=\frac{1}{2\lambda}$ and $\langle u,v \rangle=\frac{1}{2\lambda}\|u\|_{BV}$ and finally, \label{enumb} \\
\item $\|v\|_{BV}=\frac{\mu}{2\lambda}$, $\|v\|_G=\frac{1}{2\lambda}$, $\langle u,v \rangle=\frac{1}{2\lambda}\|u\|_{BV}$ and $\langle v,w \rangle=\frac{\mu}{2\lambda}\|w\|_G$. \label{enumc} \\
\end{enumerate}
Conversely, all triplet $(u,v,w)$ which fulfills (\ref{enuma}), or (\ref{enumb}), or (\ref{enumc}) is optimal for $f=u+v+w$ and their corresponding values of $\lambda$ and $\mu$.\\
\end{thm}

An example of interest of Theorem \ref{theo:meyer1} is given by the following observation. Let us assume that we have $\|f\|_G<\frac{\pi}{\lambda\mu}$. Then, the optimal decomposition is given by case (\ref{enuma}). Indeed, if we compare the optimal decomposition $f=u+v+w$ to the trivial decomposition $f=0+0+f$, we have
\begin{equation}
\|u\|_{BV}+\lambda \|v\|_{L^2}^2+\mu \|w\|_G \leqslant \mu \|f\|_G
\end{equation}
which implies
\begin{equation}
\|v\|_{L^2}\leqslant \sqrt{\frac{\mu}{\lambda}\|f\|_G}.
\end{equation}
But $\|v\|_G\leqslant \frac{1}{2\sqrt{\pi}}\|v\|_{L^2}$ which results in $\|v\|_G<\frac{1}{2\lambda}$.\\

\begin{IEEEproof}
Let us return to the proof of Theorem \ref{theo:meyer1}. First, let us observe that in the case $0<\mu<4\pi$ the problem is equivalent to minimize

\begin{equation}
\lambda \|v\|_{L^2}^2+\mu\|f-v\|_G.
\end{equation}

Then we apply Lemma \ref{lemme:meyer1} (with $E=G$ and $E^*=BV$) to $\lambda\|v\|_{L^2}^2+\mu\|f-v\|_G$.
If $\|f\|_{BV}\leqslant \frac{\mu}{2\lambda}$, the minimum is reached if $v=0$.\\

Here, we can make a partial conclusion: if $0<\mu<4\pi$ then $u=0$. Else if $\|f\|_{BV}\leqslant \frac{\mu}{2\lambda}$, then $u=v=0$. \\

Furthermore, we need to miminize $E(u,v)=\|u\|_{BV}+\lambda\|v\|_{L^2}^2+\mu\|w\|_G$ under the constraint $f=u+v+w$ (it implies $w=f-u-v$). Assume that $v$ is fixed and we seek for the minimum with respect to $u$. If we write $\sigma=u+w$ and
\begin{equation}
\interleave \sigma \interleave =\inf \{\|u\|_{BV}+\mu\|w\|_G\; ; \; \sigma=u+w\}
\end{equation}
then
\begin{equation}
\inf_{u,v}E(u,v)=\inf_{\sigma}\{\interleave \sigma \interleave + \lambda\|f-\sigma\|_{L^2}^2\}.
\end{equation}
To minimize $\interleave \sigma \interleave + \lambda\|f-\sigma\|_{L^2}^2$, we apply Lemma \ref{lemme:meyer2}. The dual norm of $\interleave .\interleave$ is
\begin{equation}
\interleave . \interleave_*=\sup\left\{\|.\|_G \; , \; \frac{1}{\mu}\|.\|_{BV}\right\}.
\end{equation}

Our next step is the following lemma.

\begin{lem}
If $\|f\|_G\leqslant\frac{1}{2\lambda}$ and $\|f\|_{BV}\leqslant\frac{\mu}{2\lambda}$, then the minimum of $E(u,v)$ is reached for $u=w=0$ and is given by $\lambda\|f\|_{L^2}^2$.\\
\end{lem}

\begin{IEEEproof}
Indeed, as $\|f\|_G\leqslant\frac{1}{2\lambda}$ and $\|f\|_{BV}\leqslant\frac{\mu}{2\lambda}$, we get 
\begin{equation}
\interleave f\interleave_*=\sup \left\{\|f\|_G;\frac{1}{\mu}\|f\|_{BV}\right\}\leqslant\frac{1}{2\lambda}.
\end{equation}
That brings back us to the case (\ref{enumaa}) of Lemma \ref{lemme:meyer1}. Then $\sigma=0$ and $v=f$ and with the agreement to the definition of the norm $\interleave . \interleave$, this implies $u=w=0$.
\end{IEEEproof}
This result concludes the first assertion in Theorem \ref{theo:meyer1}. Now, let us look to the second assertion where
\begin{equation}\label{hyp12}
\|f\|_G\leqslant\frac{1}{2\lambda} \quad \text{and} \quad \|f\|_{BV}>\frac{\mu}{2\lambda}.
\end{equation}

Observing that $\|f\|_G\leqslant\frac{1}{4\pi}\|f\|_{BV}$. Then we can't have $\|f\|_G>\frac{1}{2\lambda}$ and $\|f\|_{BV}\leqslant\frac{\mu}{2\lambda}$ if $0<\mu \leqslant 4\pi$.\\
Under the assumption (\ref{hyp12}) Lemma \ref{lemme:meyer1} ensures us that the optimal $\sigma$ fulfills
\begin{equation}
\interleave v\interleave_*=\frac{1}{2\lambda} \quad \text{and} \quad \langle v,\sigma\rangle=\frac{1}{2\lambda}\interleave \sigma\interleave.
\end{equation}
In addition $\interleave\sigma\interleave=\|u\|_{BV}+\mu\|w\|_G$ because $u$ and $w$ are optimized. We either have
\begin{equation}\label{hyp13}
\|v\|_{BV}=\frac{\mu}{2\lambda}\quad\text{and}\quad\|v\|_G<\frac{1}{2\lambda}
\end{equation}
or
\begin{equation}\label{hyp14}
\|v\|_{BV}\leqslant\frac{\mu}{2\lambda}\quad\text{and}\quad\|v\|_G=\frac{1}{2\lambda}.
\end{equation}

Let us examine the first case. We have

\begin{lem}\label{lemme:meyer3}
If (\ref{hyp12}) and (\ref{hyp13}) are fulfilled simultaneously, then the optimal decomposition $f=u+v+w$ verifies $u=0$ and $\langle v,w\rangle=\frac{\mu}{2\lambda}\|w\|_G$.\\
\end{lem}

\begin{IEEEproof}
Indeed, Lemma \ref{lemme:meyer1} yields
\begin{equation}\label{eq15}
\langle v,u+w\rangle=\frac{1}{2\lambda}\left(\|u\|_{BV}+\mu\|w\|_G\right).
\end{equation}
But
\begin{equation}
\langle v,w\rangle\leqslant \|v\|_{BV}\|w\|_G=\frac{\mu}{2\lambda}\|w\|_G
\end{equation}
while
\begin{equation}
\langle v,u\rangle\leqslant\|v\|_G\|u\|_{BV}<\frac{1}{2\lambda}\|u\|_{BV}.
\end{equation}
Adding these two inequalities, we get (\ref{eq15}). These inequalities must be equalities. 
This implies $u=0$ and $\langle v,w\rangle=\frac{\mu}{2\lambda}\|w\|_G$.
\end{IEEEproof}

Now look at the second case. We distinguish
\begin{equation}\label{eq16}
\|v\|_{BV}<\frac{\mu}{2\lambda}\quad\text{and}\quad\|v\|_G=\frac{1}{2\lambda}
\end{equation}
and
\begin{equation}\label{eq17}
\|v\|_{BV}=\frac{\mu}{2\lambda}\quad\text{and}\quad\|v\|_G=\frac{1}{2\lambda}.
\end{equation}
From (\ref{eq16}), by the same arguments previously used, we can conclude that $w=0$. Then $f=u+v$ is the optimal decomposition and we have (see Lemma \ref{lemme:meyer1}) $\|v\|_G=\frac{1}{2\lambda}$ et $\langle u,v\rangle=\|u\|_{BV}\frac{1}{2\lambda}$.\\
Let us address the converse implication. The result is given by the following lemma.

\begin{lem}\label{lem:lem7}
Let us assume that $\|v_0\|_G=\frac{1}{2\lambda}$ and $\|v_0\|_{BV}\leqslant\frac{\mu}{2\lambda}$ where $\langle u_0,v_0\rangle=\frac{1}{2\lambda}\|u_0\|_{BV}$. We write $f_0=u_0+v_0$.\\
Then for all functions $\alpha\in BV$ and all functions $w\in L^2(\R^2)$, we have
\begin{align}\label{eq18}
\|u_0+\alpha\|_{BV}+&\lambda\|v_0-\alpha-w\|_{L^2}^2+\mu\|w\|_G \\ \notag &\geqslant \|u_0\|_{BV}+\lambda\|v_0\|_{L^2}^2.
\end{align}
\end{lem}
This means that for such a function $f_0$ and for these parameters $\lambda$ and $\mu$, the decomposition $u_0+v_0$ is optimal.\\

\begin{IEEEproof}
For proving Lemma \ref{lem:lem7}, we divide (\ref{eq18}) by $2\lambda=\|v_0\|_G^{-1}$ and we get
\begin{align*}
\|u_0&+\alpha\|_{BV}\|v_0\|_G+\frac{1}{2}\|v_0-\alpha-w\|_{L^2}^2+\|v_0\|_{BV}\|w\|_G \\
=&\|u_0+\alpha\|_{BV}\|v_0\|_G+\frac{1}{2}\|v_0\|_{L^2}^2-\langle v_0,\alpha\rangle+\frac{1}{2}\|\alpha\|_ {L^2}^2 \\&-\langle w,v_0-\alpha\rangle +\frac{1}{2}\|w\|_{L^2}^2+\|v_0\|_{BV}\|w\|_G \\
\geqslant &\langle u_0,v_0\rangle+\langle \alpha,v_0\rangle+\frac{1}{2}\|v_0\|_{L^2}^2-\langle \alpha,v_0\rangle+\frac{1}{2}\|\alpha\|_{L^2}^2\\& -\langle w,v_0-\alpha\rangle+\frac{1}{2}\|w\|_{L^2}^2+\|v_0\|_{BV}\|w\|_G\\
=&\frac{1}{2\lambda}\|u_0\|_{BV}+\frac{1}{2}\|v_0\|_{L^2}^2+\frac{1}{2}\|\alpha\|_{L^2}^2-\langle w,v_0-\alpha\rangle\\&+\frac{1}{2}\|w\|_{L^2}^2+\|v_0\|_{BV}\|w\|_G\\
=&\frac{1}{2\lambda}\|u_0\|_{BV}+\frac{1}{2}\|v_0\|_{L^2}^2+\frac{1}{2}\|\alpha+w\|_{L^2}^2-\langle w,v_0\rangle\\&+\|v_0\|_{BV}\|w\|_G \\
\geqslant& \frac{1}{2\lambda}\|u_0\|_{BV}+\frac{1}{2}\|v_0\|_{L^2}^2.\\
\end{align*}

Let us notice that if we have an equality, we necessarily have $\alpha=-w$ and $\langle w,v_0\rangle=\|v_0\|_{BV}\|w\|_G$. We also get $\mu=2\lambda\|v_0\|_{BV}$. Let us return to (\ref{eq18}) which can be written
\begin{equation}
\|u_0-w\|_{BV}+\lambda\|v_0\|_{L^2}^2+\mu\|w\|_G \geqslant \|u_0\|_{BV}+\lambda\|v_0\|_{L^2}^2
\end{equation}
or
\begin{equation}
\|u_0-w\|_{BV}+2\lambda\|v_0\|_{BV}\|w\|_G \geqslant \|u_0\|_{BV}
\end{equation}
i.e
\begin{equation}
\|u_0-w\|_{BV}\|v_0\|_G+\langle w,v_0\rangle \geqslant \|u_0\|_{BV}\frac{1}{2\lambda}.
\end{equation}
But we have
\begin{equation}
\|u_0-w\|_{BV}\|v_0\|_G\geqslant \langle u_0-w,v_0\rangle=\langle u_0,v_0\rangle-\langle w,v_0\rangle.
\end{equation}
Finally
\begin{equation}
\langle u_0,v_0\rangle=\|u_0\|_{BV}\frac{1}{2\lambda}.
\end{equation}
If we have the equality, we must have
\begin{equation}
\|u_0-w\|_{BV}\|v_0\|_G=\langle u_0-w,v_0\rangle.
\end{equation}
\end{IEEEproof}

Let us examine the reciprocal of the lemma \ref{lemme:meyer3}.

\begin{lem}
Assume that $f=v_0+w_0$ with $\|v_0\|_{BV}=\frac{\mu}{2\lambda}$, $\|v_0\|_G<\frac{1}{2\lambda}$ and $\langle v_0,w_0\rangle=\frac{\mu}{2\lambda}\|w_0\|_G$.\\
Then $f=v_0+w_0$ is the optimal decomposition.\\
\end{lem}

\begin{IEEEproof}
If we write $w=w_0+\tilde{w}$ and $v=v_0$, (\ref{eq:pbmeyer}) is equivalent to
\begin{equation} 
\|u\|_{BV}+\lambda\|v_0-u-\tilde{w}\|_{L^2}^2+\mu\|w_0+\tilde{w}\|_G=J(u,\tilde{w}).
\end{equation}
Then, it comes

\begin{equation}
\|w_0+\tilde{w}\|_G\|v_0\|_{BV}\geqslant \langle w_0+\tilde{w},v_0\rangle=\langle w_0,v_0\rangle+\langle \tilde{w},v_0\rangle
\end{equation}
and, by assumption, $\langle w_0,v_0\rangle=\frac{\mu}{2\lambda}\|w_0\|_G$. Then
\begin{equation}
\|w_0+\tilde{w}\|_G\|v_0\|_{BV}\geqslant\frac{\mu}{2\lambda}\|w_0\|_G+\langle\tilde{w},v_0\rangle
\end{equation}
and as $\|v_0\|_{BV}=\frac{\mu}{2\lambda}$, we deduce that
\begin{equation}
\|w_0+\tilde{w}\|_G\geqslant\|w_0\|_G+\frac{2\lambda}{\mu}\langle \tilde{w},v_0\rangle.
\end{equation}
In addition,
\begin{align*}
\lambda\|v_0&-u-\tilde{w}\|_{L^2}^2=\\& \lambda\|v_0\|_{L^2}^2-2\lambda\langle \tilde{w},v_0\rangle-2\lambda\langle u,v_0\rangle+\lambda\|u+\tilde{w}\|_{L^2}^2\\
=&\langle v_0-u-\tilde{w},v_0-u-\tilde{w}\rangle=\langle v_0,v_0\rangle-\langle v_0,u\rangle\\&-\langle v_0,\tilde{w}\rangle-\langle u,v_0\rangle+\langle u,u\rangle+\langle u,\tilde{w}\rangle-\langle \tilde{w},v_0\rangle\\&+\langle\tilde{w},u\rangle+\langle\tilde{w},\tilde{w}\rangle\\
=&\|v_0\|_{L^2}^2-2\langle u,v_0\rangle-2\langle v_0,\tilde{w}\rangle+2\langle u,\tilde{w}\rangle+\langle u,u\rangle\\&+\langle\tilde{w},\tilde{w}\rangle.
\end{align*}
But $2\langle u,\tilde{w}\rangle+\langle u,u\rangle+\langle\tilde{w},\tilde{w}\rangle=\|u+\tilde{w}\|_{L^2}^2$ which implies
\begin{align}
J(u,\tilde{w})\geqslant&\lambda\|v_0\|_{L^2}^2+\mu\|w_0\|_G+\|u\|_{BV}-2\lambda\langle u,v_0\rangle\\ \notag
&+\lambda\|u+\tilde{w}\|_{L^2}^2.
\end{align}
To conclude, Lemma \ref{lemme:bv1} yields
\begin{equation}
\left|\langle u,v_0 \rangle\right|\leqslant \|u\|_{BV}\|v_0\|_G<\frac{1}{2\lambda}\|u\|_{BV}
\end{equation}
and
\begin{align}
\|u+w\|_{L^2}^2=\|f-v_0-&w_0\|_{L^2}^2=0 \\ \notag &(\text{we recall that}\; f=v_0+w_0).
\end{align}
\end{IEEEproof}
At this stage, we have proved the points (\ref{enuma}) and (\ref{enumb}) of Theorem \ref{theo:meyer1}.
To finish the proof, we need to establish point (\ref{enumc}). The direct part is proved by the same arguments we used in (\ref{enuma}) or (\ref{enumb}). Let us examine the reciprocal. For functions $\alpha\in BV$ and $\beta\in L^2$ choosen arbitrarily, we want to prove that $E(\alpha,\beta)\geq E(0,0)$ (see (\ref{eq:pbmeyer}) for the definition of $E(.,.)$). Then we need to calculate
\begin{equation}
\|u+\alpha\|_{BV}+\lambda\|v+\beta\|_{L^2}^2+\mu\|w-\alpha-\beta\|_G.
\end{equation}

We know that $u,v$ and $w$ verify point (\ref{enumc}).
As $\|v\|_{BV}=\frac{\mu}{2\lambda}$ and $\|v\|_{BV}\|w-\alpha-\beta\|_G\geqslant \langle v,w-\alpha-\beta\rangle$, we have
\begin{equation}
\mu\|w-\alpha-\beta\|_G\geqslant 2\lambda\left(\langle v,w\rangle-\langle v,\alpha\rangle - \langle v,\beta\rangle \right).
\end{equation}

In addition, as $\|v\|_G=\frac{1}{2\lambda}$, we also have
\begin{equation}
\|u+\alpha\|_{BV}\geqslant 2\lambda \left(\langle u,v\rangle+\langle\alpha,v\rangle\right).
\end{equation}

Finally
\begin{align}
\lambda\|v+\beta\|_{L^2}^2&=\lambda\|v\|_{L^2}^2+2\lambda\langle v,\beta\rangle+\lambda\|\beta\|_{L^2}^2, \\
\langle v,w\rangle&=\frac{\mu}{2\lambda}\|w\|_G\quad\text{and}\\
\langle v,u\rangle&=\frac{1}{2\lambda}\|u\|_{BV}.
\end{align}
This permits to conclude that all the terms disapear and that only remains (the minimum is reached for $\beta=0$) 
\begin{equation}
\|u\|_{BV}+\lambda\|v\|_{L^2}^2+\mu\|w\|_G.
\end{equation}

This ends the proof of Theorem \ref{theo:meyer1}.
\end{IEEEproof}

Let us notice that we do not have uniqueness of the decomposition. The following counterexample will persuade us.\\

Let us denote by $\theta$ the characteristic function of the unit disk, then we have $\|\theta\|_G=\frac{1}{2}$ and let us take $f=3\theta$. Let us consider the two decompositions $f=\theta+\theta+\theta$ and $f=2\theta+\theta+0$. We assume, without loss of generality, that $\lambda=1$ and $\mu=4\pi$.\\
For the first one we have $\|v\|_{BV}=\frac{\mu}{2\lambda}$, $\|v\|_G=\frac{1}{2\lambda}$, $\langle u,v\rangle=\pi=\frac{\|u\|_{BV}}{2\lambda}$ and $\langle v,w\rangle=\pi=\frac{\mu}{2\lambda}\|w\|_G$.\\
For the second one, we have well $\|v\|_{BV}\leqslant \frac{\mu}{2\lambda}$, $\langle u,v\rangle=\frac{1}{2\lambda}\|u\|_{BV}$ and $\|v\|_G=\frac{1}{2\lambda}$. We conclude that these two decompositions respect Theorem \ref{theo:meyer1} and then we don't have the uniqueness of the decomposition.\\

\section{Application}\label{sec:app}
In this section, we present an application of Theorem \ref{theo:meyer1}. Let us assume that we deal with long and thin objects in an image. This kind of object can be modeled by
\begin{equation}
f(x_1,x_2)=1 \qquad \text{if} \qquad 0 \leqslant x_1 \leqslant L , 0 \leqslant x_2 \leqslant \epsilon
\end{equation}
where $L \gg 1$ and $0<\epsilon \ll 1$. Then $\|f\|_G \leqslant \epsilon$ while $\|f\|_{BV}=2(L+\epsilon)$. It's easy to see that we are in the case \ref{enuma} of Theorem \ref{theo:meyer1} if $\epsilon<\sqrt{\frac{\pi}{\lambda\mu}}$ and $\mu<4\lambda(L+\epsilon)$ i.e if $L$ is rather large compared to $\mu$.\\

By Theorem \ref{theo:meyer1}, we conclude that $u=0$, $\|v\|_{BV}=\frac{\mu}{2\lambda}$. Then we have $\|w\|_{BV}\geqslant\|f\|_{BV}-\|v\|_{BV}\geqslant 2(L+\epsilon)-\frac{\mu}{2\lambda}$ which is high. In this case, the $w$ part is the most important one. This means that this kind of objects will be attracted in the $w$ component.\\

\begin{figure}[!t]
\begin{center}
\begin{tabular}{cc}
\includegraphics[width=0.47\textwidth]{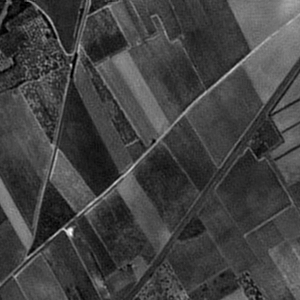} \\
\includegraphics[width=0.47\textwidth]{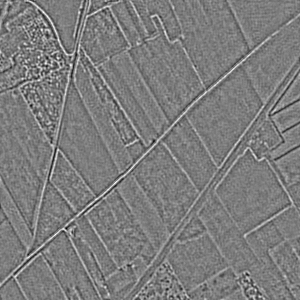}
\end{tabular}
\end{center}
\caption{Example of a portion of an aerial image: original image on top, $w$ component on bottom which lets appear enhanced roads.}
\label{fig:rehaussement}
\end{figure}

\begin{figure}[!t]
\begin{center}
\includegraphics[width=0.47\textwidth]{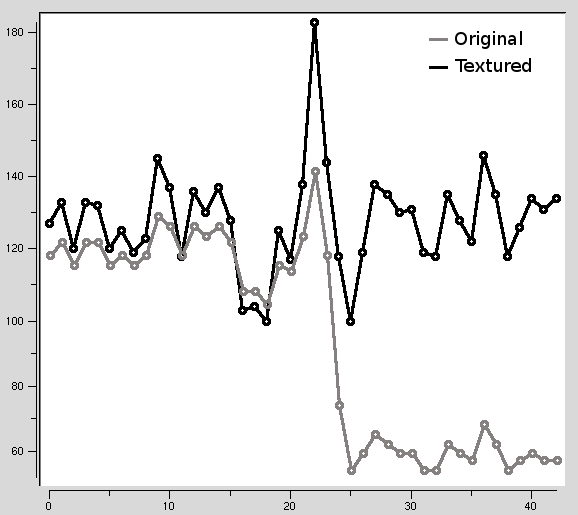}
\end{center}
\caption{Illustration of the constrast enhancement of thin and long structures.}
\label{fig:profil}
\end{figure}

\begin{figure}[!t]
\begin{center}
\begin{tabular}{cc}
\includegraphics[width=0.47\textwidth]{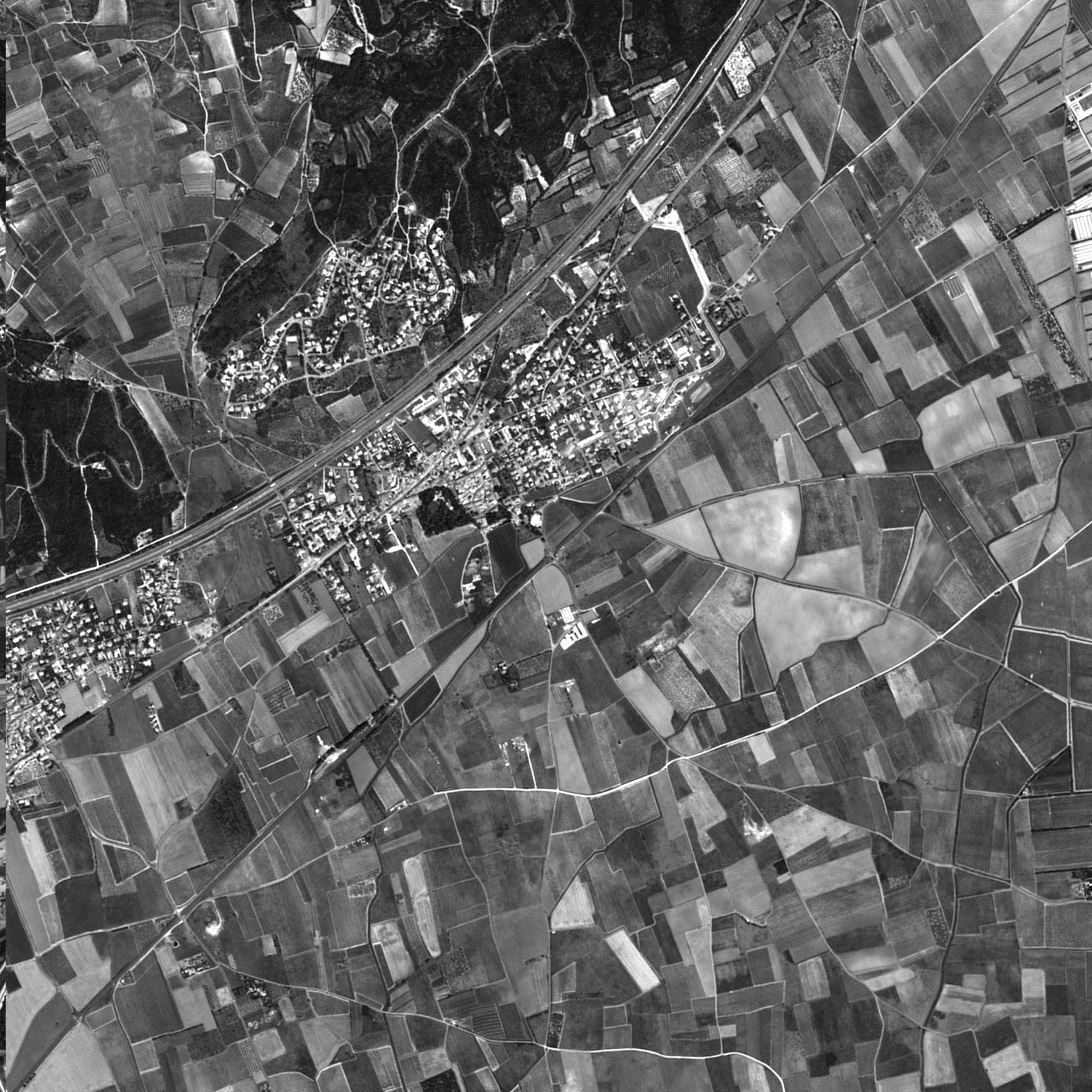} \\
\includegraphics[width=0.47\textwidth]{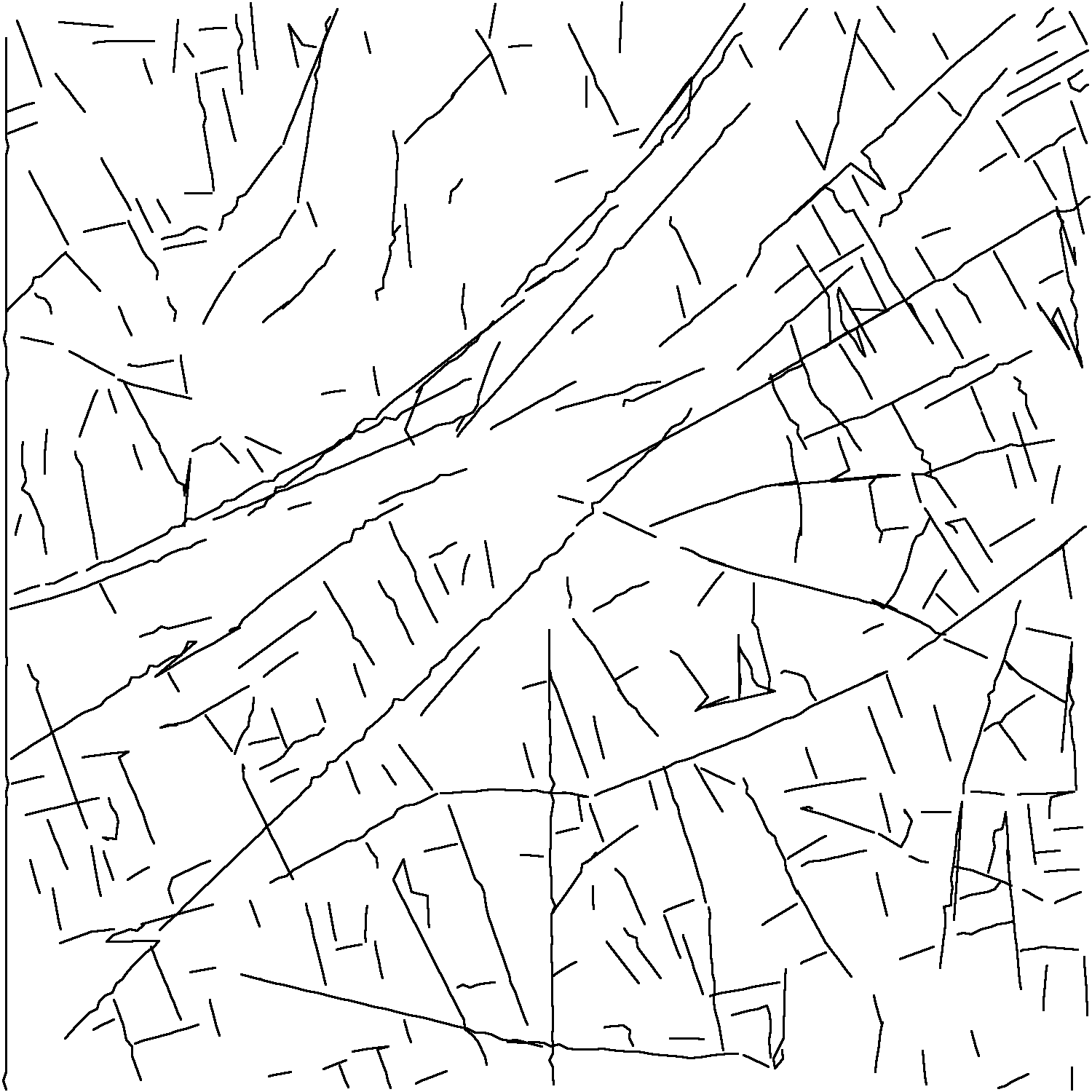}
\end{tabular}
\end{center}
\caption{Example of road network detection: original image on top, detected roads on bottom.}
\label{fig:road}
\end{figure}

This property was used in \cite{jegilles2} as a preprocessing stage in an aerial road networks detection. Indeed, road networks could be considered as long and thin objects in the image. The previous result teaches us that this kind of objects will be enhanced in the texture component (but $u$ is not strictly equal to $0$ and $w$ does not contain only roads because the original image contains different kind of objects). So we decompose the image and then apply a detection algorithm on the $w$ component. Figure \ref{fig:rehaussement} shows a zoomed portion of an aerial image and its $w$ component, figure \ref{fig:profil} exhibits the same contrast evolution from one side to another of a road in the original image and the texture component respectively. We clearly see that roads are the most visible objects in the texture component. Figure \ref{fig:road} shows a result we get, on a bigger image, by this approach with a very simple detection algorithm applied on the $w$ component (see the appendix for details about the practical algorithm used).

\section{Conclusion}
In this paper, we present some new theoretical results about the second author's $BV-G$ decomposition model which separates structures and textures from an original image. The main theorem we proved gives the optimal decomposition we get with regard to the parameter selection and the properties of the images, in terms of the different norms associated to the involved function spaces. This theorem permits us to propose an enhancement method for long and thin objects before a detection stage. This method was tested on an aerial or satellite road networks detection application.\\

In a future work, the results of the main theorem could be associated with different kind of objects and could permit to help in selecting values of the algorithm's parameters in order to get good decomposition results. Another way of research is to extend these results to other function spaces (like Besov spaces) or to three parts decomposition models which deal with noise (as cited in Section \ref{sec:tvg}).

\appendix
\section{Practical algorithms}
In this appendix we briefly recall the practical algorithms used in section \ref{sec:app}. The whole algorithm can be split into three parts. First, the decomposition which provides us the $w$ component on which the detection is made. The second part consists on a first stage of detection based on segments detection. The last part is a refinement stage which converts the previously detected segments into active contours in order to get the real topology of roads. Let us give more details on each part.\\

The image decomposition algorithm is the one proposed proposed by Aujol in \cite{aujol,aujolphd} based on Chambolle's projectors cited in section \ref{sec:tvg}. The slightly modified model of Aujol is defined by equation (\ref{eq:aujolg}) (Aujol proved that the minimizers of its model are also minimizers of the original model of the second author). 

\begin{equation}\label{eq:aujolg}
F_{\lambda,\mu}^{AU}(u,v)=J(u)+J^*\left(\frac{v}{\mu}\right)+(2\lambda)^{-1}\|f-u-v\|_{L^2}^2
\end{equation}

where
\begin{equation}
(u,v)\in BV(\Omega)\times G_{\mu}(\Omega).
\end{equation}

and the set $G_{\mu}$ is the subset of $G$ where $\forall v\in G_{\mu}$, $\|v\|_G\leqslant\mu$. Moreover, $J^*$ is the characteristic function over $G_1$ with the property that $J^*$ is the dual operator of $J$ ($J^{**}=J$). Thus,

\begin{equation}\label{eq:jdual}
J^*(v)=
\begin{cases}
0 \qquad \text{if} \; v\in G_1 \\
+\infty \qquad \text{else.}
\end{cases}
\end{equation}

The minimizers can be found by the following iterative algorithm.

\begin{center}
\begin{enumerate}
\item Initialization:
\begin{equation*}
u_0=v_0=0
\end{equation*}
\item Iteration $n+1$:
\begin{eqnarray*}
&v_{n+1}=P_{G_{\mu}}(f-u_n)\\
&u_{n+1}=f-v_{n+1}-P_{G_{\lambda}}(f-v_{n+1})
\end{eqnarray*}
\item We stop the algorithm if
\begin{equation*}
\max\left(|u_{n+1}-u_n|,|v_{n+1}-v_n|\right)\leqslant \epsilon
\end{equation*} 
or if we reach a prescribed maximal number of iterations.\\
\end{enumerate}
\end{center}

The expressions of Chambolle's projectors can be found in \cite{aujol,aujolphd} and are very easy to implement.\\

The first stage of the detection algorithm is the one proposed by Morel's team in \cite{gestalt1}. It is based on an a contrario formulation issued from the Gestalt theory. The output of this algorithm is a set of segments corresponding to aligned points in the image.\\

\begin{figure}[!t]
\includegraphics[width=0.47\textwidth]{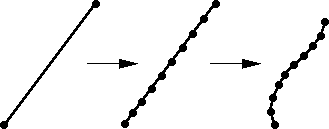}
\caption{Segment to polygonal active contour conversion strategy.}
\label{fig:snake}
\end{figure}

In the last stage, we start by filtering the set of segments. We mean that we fusion very close segments, we supplement each segment which follows another one. Then each segment is converted into an open polygonal active contour (see figure \ref{fig:snake}). As they are very close to the final position (we recall that it is a refinement stage), we can use the active contour algorithm proposed by the first author in \cite{jegilles3}.

\begin{thebibliography}{10}
\bibitem{aujol}
{\sc J.F. Aujol, G. Aubert, L. Blanc-F\'eraud and A. Chambolle},
{\em Image decomposition into a bounded variation component and an oscillating component},
Journal of Mathematical Imaging and Vision, vol.22, No.1, pp.~71--88, 2005.

\bibitem{aujolphd}
{\sc J.F.Aujol},
{\em Contribution \`a l'analyse de textures en traitement d'images par m\'ethodes variationnelles et \'equations aux d\'eriv\'ees partielles},
Nice Sophia Antipolis University, Ph.D Thesis, 2004.

\bibitem{aujol2}
{\sc J.F.Aujol, G.Aubert, L.Blanc-F{\'e}raud and A. Chambolle},
{\em Decomposing an image: Application to {SAR} images},
Scale-Space '03, Lecture Notes in Computer Science, vol.1682, 2003.
 
\bibitem{vese1}
{\sc L.A.Vese and S.J.Osher},
{\em Modeling textures with total variation minimization and oscillating patterns in image processing},
Journal of Scientific Computing, vol.19, pp.~553--572, 2003.

\bibitem{vese2}
{\sc S.J.Osher, A.Sole and L.A.Vese},
{\em Image decomposition and restoration using total variation minimization and the {H}$^{-1}$ norm},
Multiscale Modeling and Simulation: A SIAM Interdisciplinary Journal, vol.1, No.3, pp.~349--370, 2003.

\bibitem{meyer}
{\sc Yves Meyer},
{\em Oscillating patterns in image processing and in some nonlinear evolution equations},
The Fifteenth Dean Jacquelines B. Lewis Memorial Lectures, American Mathematical Society, March 2001.

\bibitem{chambolle}
{\sc A. Chambolle},
{\em An algorithm for total variation minimization and applications},
Journal of Mathematical Imaging and Vision, vol.20, No.1-2, pp.~89--97, 2004.

\bibitem{rof}
{\sc L.Rudin, S.Osher and E.Fatemi},
{\em Nonlinear total variation based noise removal algorithms},
Physica D, vol.60, pp.~259--268, 1992.

\bibitem{aujoluvw}
{\sc J.F.Aujol and A.Chambolle},
{\em Dual norms and image decomposition models},
International Journal of Computer Vision, vol.63, No.1, pp.~85--104, 2005.

\bibitem{jegilles}
{\sc J.Gilles},
{\em Noisy image decomposition: a new structure, textures and noise model},
Journal of Mathematical Imaging and Vision, vol.28, No.3, pp.~285--295, 2007.

\bibitem{jegilles2} 
{\sc J.Gilles},
{\em Choix d'un espace de repr\'esentation image adapt\'e \`a la d\'etection de r\'eseaux routiers},
TAIMA (Traitement et Analyse de l'Information: M\'ethodes et Applications) Conference, Hammamet, Tunisia, 2007.

\bibitem{aujolcouleur}
{\sc J.F.Aujol and S.H.Kang},
{\em Color image decomposition and restoration},
Journal of Visual Communication and Image Representation, vol.17, No.4, pp.~916--928, 2006.

\bibitem{aujolclassiftexture}
{\sc J.F.Aujol and T.Chan},
{\em Combining geometrical and textured information to perform image classification},
Journal of Visual Communication and Image Representation, vol.17, No.5, pp.~1004--1023, 2006.

\bibitem{triet1}
{\sc Triet M. Le and Luminita A. Vese},
{\em Image Decomposition Using Total Variation and div(BMO)},
Multiscale Modeling and Simulation: A SIAM Interdisciplinary Journal, vol.4, No.2, pp.~390--423, 2005


\bibitem{triet2}
{\sc John B. Garnett, Triet M. Le, Yves Meyer, Luminita A. Vese},
{\em Image decompositions using bounded variation and generalized homogeneous Besov spaces},
Appl. Comput. Harmon. Anal., vol.23, pp.~25--56, 2007

\bibitem{gestalt1}
{\sc A. Desolneux, L. Moisan and J.-M. Morel},
{\em Maximal Meaningful Events and Applications to Image Analysis},
Annals of Statistics, vol.31, No.6, pp.~1822--1851, 2003

\bibitem{jegilles3}
{\sc J. Gilles and B. Collin},
{\em Fast probabilist snake algorithm},
in Proceedings of International Conference on Image Processing (ICIP), Barcelona, 2003

\bibitem{chambolle2}
{\sc A. Chambolle, V. Caselles, M. Novaga, D. Cremers and T. Pock},
{\em An introduction to total variation for image analysis},
in Lecture Notes for the summer school on sparsity in Linz, Austria, August 2009, available at http://hal.archives-ouvertes.fr/docs/00/43/75/81/PDF/preprint.pdf

\end{thebibliography}
\end{document}